\documentstyle[amssymb,amsfonts,12pt]{amsart}

\newtheorem{theorem}{Theorem}[section]

\newtheorem{proposition}[theorem]{Proposition}

\newtheorem{question}[theorem]{Question}

\theoremstyle{remark}

\def\QSet{\mbox{\rm\kern.24em
\vrule width.03em height1.48ex depth-.051ex \kern-.26em Q}}

\def\x{{\bf x}}\def\y{{\bf y}}\def\z{{\bf z}}\def\0{{\bf 0}}

\def\Z{{\mathbb Z}}

\def\L{{\mathcal L}}
\def\det{{\operatorname{det}}}

\def\bas{\begin{align*}}
\def\eas{\end{align*}}
\def\bi{\begin{itemize}}
\def\ei{\end{itemize}}

\def \endprf{\hfill  {\vrule height6pt width6pt depth0pt}\medskip}
\def\emph#1{{\it #1}}

\begin{document}

\title[On the number of non-congruent lattice tetrahedra ]{On the number of non-congruent lattice tetrahedra}

\author{Ciprian Demeter}
\address{Department of Mathematics, Indiana University, 831 East 3rd St., Bloomington IN 47405}
\email{demeterc@@indiana.edu}

\keywords{}
\thanks{The author is supported by  the NSF Grant DMS-1161752}

\begin{abstract}
We prove that there are  "many" non-congruent tetrahedra in the truncated lattice $[0,q]^3\cap \Z^3$. This answers a question from \cite{GILP}.
\end{abstract}
\maketitle

\section{Introduction}

Our goal is to answer the following question asked in the last section of \cite{GILP}. 
\begin{question}
Let $T_3([0,q]^3\cap \Z^3)$ denote the collection of all equivalence classes of congruent tetrahedra with vertices in $[0,q]^3\cap \Z^3$. Is there a  $\delta>0$ and some $C>0$, both independent of $q$ such that
$$\#T_3([0,q]^3\cap \Z^3)\le Cq^{9-\delta}$$
for each $q>1$?
\end{question}
A positive answer to this question would have implications on producing lower bounds for the Falconer distance-type problem for tetrahedra. We refer to \cite{GILP} for details on this interesting problem.

Here we give a negative answer to this question. We hope that our approach to answering this question will inspire further constructions which might eventually improve the lower bound.
\begin{theorem}
\label{mainthmmmm}
We have for each $\epsilon>0$ and each $q>1$
$$\#T_3([0,q]^3\cap \Z^3)\gtrsim_{\epsilon} q^{9-\epsilon}.$$
\end{theorem}
Note the following trivial upper bound, which shows the essential tightness of our result
$$\#T_3([0,q]^3\cap \Z^3)\le Cq^{9}.$$
Indeed,  by translation invariance it suffices to fix one vertex at the origin. The upper bound follows since there are $(q+1)^3$ possibilities for each of the remaining three vertices.   

\section{Systems of quadratic equations}
The proof of Theorem \ref{mainthmmmm} will rely on the theory developed in \cite{BD} around Siegel's mass formula. To keep the presentation essentially self-contained, we briefly recall the relevant theorems here and refer the interested reader to \cite{BD} for further details.

Let $m=n+1\ge 3$ and let $\gamma\in M_{m,m}(\Z)$ and $\Lambda\in M_{n,n}(\Z)$ be two positive definite matrices with integer entries. Denote by $A(\gamma,\Lambda)$  the number of  solutions $\L\in M_{m,n}(\Z)$ for
\begin{equation}
\label{Snew27}
\L^*\gamma\L=\Lambda.
\end{equation}

For primes $p$ define
$$\nu_p(\gamma,\Lambda)=\lim_{r\to\infty}\frac1{p^{r(mn-\frac{n(n+1)}{2})}}|\{\L\in M_{m,n}(\Z_{p^r}):\;\L^*\gamma\L\equiv \Lambda \mod p^r\}|.$$
The following is an immediate consequence of the Siegel's mass formula in \cite{Si}.
\begin{equation}
\label{Snew1863}
A(\gamma,\Lambda)\lesssim_{n,\gamma}\prod_{p\text{ prime}}\nu_p(\gamma,\Lambda).
\end{equation}
In our  application here $m=3, n=2$ and $\gamma$ is the identity matrix $I_{3}$.
\bigskip

Fix $\Lambda\in M_{n,n}(\Z)$, a  positive definite matrix, in particular $\det(\Lambda)\not=0$. In evaluating $\nu_p(I_{n+1},\Lambda)$ we distinguish two separate cases: $p\nmid\det(\Lambda)$ and $p|\det(\Lambda)$.

\begin{proposition}
\label{propnondiv}
Assume $p$ is not a factor of $\det (\Lambda)$. Then
$$\nu_p(I_{n+1},\Lambda)\le 1+\frac{C}{p^2},$$
where $C$ is independent of $p,\Lambda.$
\end{proposition}

The contribution of these primes to the product \eqref{Snew1863} is easily seen to be $O(1)$.

For an $n\times n$ matrix $\Lambda$ and for $A,B\subset \{1,\ldots,n\}$ with $|A|=|B|$ we define
$$\mu_{A,B}=\det((\Lambda_{i,j})_{i\in A,j\in B}).$$
Also, for an integer $k$ and a prime number $p$ we denote by $o_p(k)$ the largest positive integer $\alpha$ such that $p^\alpha|k$.

\begin{proposition}
\label{thecaseofpdivisor}
Assume $p|\det(\Lambda)$. Then there is $C$ independent of $p, \Lambda$ such that
$$\nu_p(I_{n+1},\Lambda)\le C \sum_{0\le l_i:1\le i\le n\atop{l_1+l_2+\ldots+l_n\le o_p(\det(\Lambda))}}p^{\beta_2(l_1,\ldots,l_n)+\ldots+\beta_n(l_1,\ldots,l_n)},$$
where $\beta_i=\beta_i(l_1,\ldots,l_n)$ satisfies
$$\beta_i=\min\{(i-1)l_i,(i-2)l_i+\min_{|A|=1}o_p(\mu_{\{1\},A})-l_1,(i-3)l_i+\min_{|A|=2}o_p(\mu_{\{1,2\},A})-l_1-l_2,\ldots,$$$$\ldots,\min_{|A|=i-1}o_p(\mu_{\{1,2,\ldots,i-1\},A})-l_1-l_2-\ldots-l_{i-1}\}$$
\end{proposition}

\section{Proof of Theorem \ref{mainthmmmm}}
In the following discussion, unless specified otherwise, $\epsilon$ will denote an arbitrary positive number.

It will be useful to  recall the classical facts that  circles of radius $r$ contain $O(r^\epsilon)$ lattice points while spheres  of radius $r$ centered at the origin contain $O(r^{1+\epsilon})$ lattice points.

As mentioned before, we fix one vertex to be the origin $\0=(0,0,0)$. A class of congruent tetrahedra in $T_3([0,q]^3\cap \Z^3)$ can be identified with a matrix $\Lambda\in M_{3,3}(\Z)$. Namely, the congruence class of the tetrahedron with vertices $\0,\x,\y,\z\in [0,q]^3\cap \Z^3$ is represented by the matrix 
$$\Lambda=\begin{bmatrix}\langle \x,\x\rangle &\langle \x,\y\rangle & \langle \x,\z\rangle\\ \langle \y,\x\rangle &\langle \y,\y\rangle & \langle \y,\z \rangle\\ \langle \z,\x\rangle &\langle \z,\y\rangle & \langle \z,\z\rangle\end{bmatrix}.$$
A tetrahedron is called non-degenerate if $\x,\y,\z$ are linearly independent. We will implicitly assume the congruence classes correspond to non-degenerate tetrahedra.

We seek for an upper bound on the number $N_{\Lambda}$ of integral solutions $\L=(\x,\y,\z)\in (\Z^3)^3$ to the equation
$$\L^*\L=\Lambda.$$
This will represent the number of congruent tetrahedra with  side lengths specified by $\Lambda$. 
  
In the numerology from the previous section,  this corresponds to $n=m=3$. To make the theorems in that section applicable we reduce the counting problem to the $m=3,n=2$ case as follows. One can certainly bound $N_\Lambda$ by $q^{\epsilon} N'_\Lambda$,  where $N'_\Lambda$ is the number of integral solutions $\L'=(\x,\y)\in (\Z^3)^2$ satisfying $$(\L')^*\L'=\Lambda'$$
and $\Lambda'$ is the $2\times 2$ minor of $\Lambda$ obtained from the first two rows and columns of $\Lambda$.
Indeed, if $\x,\y$ are fixed, the matrix $\Lambda$ forces $\z$ to lie on the intersection of the sphere of radius $\Lambda_{3,3}^{1/2}$ centered at the origin with, say, a sphere centered at $\x$ whose radius is determined only by the entries of $\Lambda$. These radii are $O(q)$, so the resulting circle can only have $O(q^\epsilon)$ points.

Note also that we only care about those $\Lambda'$ for which there exist $\x,\y\in [0,q]^3\cap \Z^3$ linearly independent, such that
$$\Lambda'=\begin{bmatrix}\langle \x,\x\rangle &\langle \x,\y\rangle\\ \langle \y,\x\rangle &\langle \y,\y\rangle\end{bmatrix}.$$
This in particular forces $\Lambda'$ to be positive definite.

Apply now Propositions \ref{propnondiv} and  \ref{thecaseofpdivisor} combined with \eqref{Snew1863} to the matrix $\Lambda'$ ($n=2$). For Proposition  \ref{thecaseofpdivisor} use $\beta_2\le \min_{|A|=1}o_p(\mu_{\{1\},A})-l_1$, then $\beta_2\le \min_{|A|=1}o_p(\mu_{\{2\},A})-l_1$ and the fact that   there are
\begin{equation}
\label{nottoomanyefgryft77856t785687}
O(\frac{\log \det(\Lambda)}{\log\log \det(\Lambda)})
\end{equation}
primes $p$ which divide  $\det(\Lambda)$.
This will bound $N'_\Lambda$ by
$$q^\epsilon \operatorname{gcd}(\Lambda_{i,j}:i,j\not=3)\le q^\epsilon \operatorname{gcd}(\Lambda_{1,1},\Lambda_{2,2}).$$
Thus
\begin{equation}
\label{dfjkgy58g0-riu0g5690i}
N_\Lambda\lesssim_\epsilon q^\epsilon \operatorname{gcd}(\Lambda_{1,1},\Lambda_{2,2}),
\end{equation}
for each $\Lambda$ corresponding to a non-degenerate tetrahedron.

Denote by $M_r$  the number of lattice points on the sphere or radius $r^{1/2}$ centered at the origin. In our case $r\le q^2$ so we know that $M_r\lesssim_\epsilon q^{1+\epsilon}$. We need to work with spheres that contain many points. Let
$$A:=\{r\le q^2:M_r\ge q/2\}.$$
Since for each $\epsilon>0$ we have $M_r\le C_\epsilon q^{1+\epsilon}$, a double counting argument shows that  $q^3\le  C_\epsilon\#Aq^{1+\epsilon}+\frac12q^2q$. Thus $\#A\gtrsim_\epsilon q^{2-\epsilon}$.

Note that for each $r_i\in A$ there are $\sim M_{r_1}M_{r_2}M_{r_3}$ non-degenerate tetrahedrons with vertices $\x,\y,\z$ on the spheres centered at the origin and with radii $r_1^{1/2},r_2^{1/2},r_3^{1/2}$ respectively. The congruence class of such a tetrahedron contains
$$\lesssim_\epsilon q^\epsilon \operatorname{gcd}(r_1,r_2)$$
elements, according to \eqref{dfjkgy58g0-riu0g5690i}.

We conclude that there are at least 
$$\frac{ M_{r_1}M_{r_2}M_{r_3}}{q^{\epsilon}\operatorname{gcd}(r_1,r_2)}$$
congruence classes generated by such non-degenerate tetrahedra. As distinct radii necessarily give rise to distinct congruence classes, we obtain the lower  bound 
$$\#T_3([0,q]^3\cap \Z^3)\gtrsim_\epsilon\sum_{r_i\in A}\frac{ M_{r_1}M_{r_2}M_{r_3}}{q^{\epsilon}\operatorname{gcd}(r_1,r_2)}\gtrsim_\epsilon q^{3-\epsilon}\sum_{r_1,r_2,r_3\in A}\frac1{\operatorname{gcd}(r_1,r_2)}.$$

It is immediate that for each integer $d$ there can be  at most $\frac{q^6}{d^2}$ tuples $(r_1,r_2,r_3)\in [0,q^2]^3$, hence also in $A^3$, with $\operatorname{gcd}(r_1,r_2)=d$. Using this observation and the bound $\#(A^3)\ge C_\epsilon q^{6-\epsilon}$, it follows that for each $\epsilon>0$ at least $\frac12\#(A^3)$ among the triples $(r_1,r_2,r_3)\in A^3$ will have  $\operatorname{gcd}(r_1,r_2)\le \frac{10q^\epsilon}{C_\epsilon}$.

This implies that
$$\sum_{r_1,r_2,r_3\in A}\frac1{\operatorname{gcd}(r_1,r_2)}\gtrsim_\epsilon q^{6-\epsilon},$$
which finishes the proof of the theorem.


\begin{thebibliography}{99}
\bibitem{BD} J. Bourgain and C. Demeter, {\em New bounds for the  discrete Fourier restriction to the sphere in four and five dimensions}, available on arxiv http://arxiv.org/pdf/1310.5244.pdf
\bibitem{GILP} A. Greenleaf, A. Iosevich, B. Liu and E. Palsson, {\em A group-theoretic viewpoint on Erd\"os-Falconer problems
and the Mattila integral}, available on arxiv  http://arxiv.org/pdf/1306.3598.pdf
\bibitem{Si} L. C. Siegel {\em Lectures on the analytical theory of quadratic forms: second term 1934-35}, The Institute for Advanced Study and Princeton University,  revised edition 1949, reprinted January 1955.
\end{thebibliography}
\end{document}